\documentclass{article}%
\usepackage{amssymb}
\usepackage{amscd}
\usepackage{graphicx}
\usepackage{amsmath}
\usepackage{amsfonts}%
\newtheorem{theorem}{Theorem}

\newtheorem{definition}[theorem]{Definition}

\newtheorem{lemma}[theorem]{Lemma}

\newtheorem{proposition}[theorem]{Proposition}
\newtheorem{remark}[theorem]{Remark}

\newenvironment{proof}[1][Proof]{\textbf{#1.} }{\ \rule{0.5em}{0.5em}}
\begin{document}

\title{The Affine Bundle Theorem in Synthetic Differential Geometry of Jet Bundles}
\author{Hirokazu Nishimura\\Institute of Mathematics, University of Tsukuba, Tsukuba\\Ibaraki 305-8571, Japan}

%change: the following command removes the date
\date{}

\maketitle

\begin{abstract}
In this paper we will establish the affine bundle theorem in our synthetic
approach to jet bundles in terms of infinitesimal spaces $D_{n}$'s, We will
then compare these affine bundles with the corresponding ones constructed in
our previous synthetic approach to jet bundles in terms of infinitesimal
spaces $D^{n}$'s both in the general microlinear setting and in the
finite-dimensional setting.

\end{abstract}

\section{Introduction}

In our \cite{n4} we have established the affine bundle theorem in our
synthetic approach to jet bundles in terms of infinitesimal spaces $D^{n}$'s.
In our succeeding \cite{n5} we have introduced another synthetic approach to
jet bundles in terms of infinitesimal spaces $D_{n}$'s, and have compared it
with the former approach both in the general microlinear setting and in the
finite-dimensional setting. However our comparison in the finite-dimensional
setting was incomplete, for our use of dimension counting techniques has
tacitly assumed that jet bundles in terms of $D^{n}$'s and $D_{n}$'s are
vector bundles, which is not the case unless the given bundle is already a
vector bundle. The principal objective in this paper is to establish the
affine bundle theorem in our second synthetic approach to jet bundles and then
to compare not merely the jet bundles based on both approaches but the affine
bundles constructed as a whole both in the general microlinear setting and in
the finite-dimensional setting. This completes our comparsion between the two
approaches in the finite-dimensional setting. Last but not least, we should
say that the notion of a simple polynomial and correspondingly the notion of a
$D^{n}$-tangential are modified essentially.

\section{Preliminaries}

\subsection{Convention}

Throughout the rest of the paper, unless stated to the contrary, $E$ and $M$
denote microlinear spaces, and $\pi:E\rightarrow M$ denotes a bundle, i.e., a
mapping from \ $E$ to $M$. The fiber of $\pi$ over $x\in M$, namely, the set
$\{y\in E\mid\pi(y)=x\}$ is denoted by $E_{x}$, as is usual. We denote by
$\mathbb{R}$ the extended set of real numbers awash in nilpotent
infinitesimals, which is expected to acquiesce in the so-called general Kock
axiom (cf. \cite{l1}). The bundle $\pi:E\rightarrow M$ is called a
\textit{vector bundle} if the fiber $E_{x}$ of $\pi$ over $x$ is an
$\mathbb{R}$-module for every $x\in M$. Various canonical injections among
infinitesimals spaces are simply denoted by the same symbol $i$. Mappings
denoted by symbols with subscripts are sometimes denoted without subscripts,
provided that the intended subscripts are clear from the context.

\subsection{Infinitesimal Spaces}

%change: the following paragraph has been embedded in a sloppypar 
%environment to prevent an overfull box}
\begin{sloppypar}
As is depicted in \cite{l1}, an \textit{infinitesimal space }corresponds to a
Weil algebra $\mathbb{R[}X_{1},...,X_{n}]/I$ in finite presentation. Given two
infinitesimal spaces $\mathbb{E}_{1}$ and $\mathbb{E}_{2}$ corresponding to
Weil algebras $\mathbb{R[}X_{1},...,X_{n}]/I$ and $\mathbb{R[}Y_{1}%
,...,Y_{m}]/J$ respectively, the infinitesimal space corresponding to the Weil
algebra $\mathbb{R[}X_{1},...,X_{n},Y_{1},...,Y_{m}]/K$, where $K$ is the
ideal generated by $I$, $J$ and $\{X_{i}Y_{j}\mid1\leq i\leq n$, $1\leq j\leq
m\}$, is denoted by $\mathbb{E}_{1}\oplus\mathbb{E}_{2}$, from which there are
canonical projections onto $\mathbb{E}_{1}$ and $\mathbb{E}_{2}$ denoted by
$p_{\mathbb{E}_{1}}$ and $p_{\mathbb{E}_{2}}$. Given two mappings
$f:\mathbb{E}_{1}\rightarrow\mathbb{F}_{1}$ and $g:\mathbb{E}_{2}%
\rightarrow\mathbb{F}_{2}$ of infinitesimal spaces, there exists a unique
mapping $f\oplus g:\mathbb{E}_{1}\oplus\mathbb{E}_{2}\rightarrow\mathbb{F}%
_{1}\oplus\mathbb{F}_{2}$ making the following diagram 
commutative:\end{sloppypar}
%change: improved diagram                                                                   
 \[ \begin{CD}
{\mathbb E}_{1}@>f>> {\mathbb F}_{1}\\
@A{p_{\mathbb E _{1}}}AA @AAp_{{\mathbb F}_{1}}A\\
{\mathbb E}_{1}\oplus {\mathbb E}_{2} @>{f\oplus g}>> {\mathbb F}_{1}
\oplus {\mathbb F}_{2}\\
@V{p_{{\mathbb E}_{2}}}VV @VV{p_{{\mathbb F}_{2}}}V\\
{\mathbb E}_{2}@>>g> {\mathbb F}_{2}
\end{CD} \]

We denote by $D_{1}$ or $D$ the totality of elements of $\mathbb{R\ }$whose
squares vanish. More generally, given a natural number $n$, we denote by
$D_{n}$ the set
\[
\{d\in\mathbb{R}|d^{n+1}=0\}\text{.}%
\]
Given natural numbers $m,n$, we denote by $D(m)_{n}$ the set
\[
\{(d_{1},...,d_{m})\in D^{m}|d_{i_{1}}...d_{i_{n+1}}=0\}\text{,}%
\]
where $i_{1},...,i_{n+1}$ shall range over natural numbers between $1$ and $m$
including both ends. We will often write $D(m)$ for $D(m)_{1}$. By convention
$D^{0}=D_{0}=\{0\}$. A polynomial $\rho$ of $d\in D_{n}$ is called a
\textit{simple} polynomial of $d\in D_{n}$ if the constant term is $0$. A
simple polynomial $\rho$ of $d\in D_{n}$ is said to be \textit{of dimension}
$m$, in notation $\mathrm{\dim}_{n}\rho=m$, provided that $m$ is the least
integer with $\rho^{m+1}=0$. By way of example, letting $d\in D_{3}$, we have
$\mathrm{\dim}_{3}d=\mathrm{\dim}_{3}(d+d^{2})=\mathrm{\dim}_{3}(d+d^{3})=3$
and $\mathrm{\dim}_{3}d^{2}=\mathrm{\dim}_{3}d^{3}=\mathrm{\dim}_{3}%
(d^{2}+d^{3})=1$. Letting $e\in D_{1}$ or $e\in D_{2}$, we have $\mathrm{\dim
}_{3}ed=1$ or $\mathrm{\dim}_{3}ed=2$ respectively. The reader should note
that our present definition of a simple polynomial is different from that in
\cite{n5}.

\textit{Simplicial infinitesimal spaces} are spaces of the form
\[
D(m;\mathcal{S})=\{(d_{1},...,d_{m})\in D^{m}|d_{i_{1}}...d_{i_{k}}=0\text{
for any }(i_{1},...,i_{k})\in\mathcal{S}\}\text{,}%
\]
where $\mathcal{S}$ is a finite set of sequences $(i_{1},...,i_{k})$ of
natural numbers with $1\leq i_{1}<...<i_{k}\leq m$. To give examples of
simplicial infinitesimal spaces, we have $D(2)=D(2;(1,2))$ and
$D(3)=D(3;(1,2),(1,3),(2,3))$, which are all symmetric. The number $m$ is
called the \textit{degree} of $D(m;\mathcal{S})$, in notation: $m=\mathrm{\deg
}D(m;\mathcal{S})$, while the maximum number $n$ such that there exists a
sequence $(i_{1},...,i_{n})$ of natural numbers of length $n$ with $1\leq
i_{1}<...<i_{n}\leq m$ containing no subsequence in $\mathcal{S}$ is called
the \textit{dimension} of $D(m;\mathcal{S})$, in notation: $n=\mathrm{\dim
}D(m;\mathcal{S})$. By way of example, $\mathrm{\deg}D(3)=\mathrm{\deg
}D(3;(1,2))=\mathrm{\deg}D(3;(1,2),(1,3))=\mathrm{\deg}D^{3}=3$, while
$\mathrm{\dim}D(3)=1$, $\mathrm{\dim}D(3;(1,2))=\mathrm{\dim}%
D(3;(1,2),(1,3))=2$ and $\mathrm{\dim}D^{3}=3$. It is easy to see that if
$n=\mathrm{\dim}D(m;\mathcal{S})$, then $d_{1}+...+d_{m}\in D_{n}$ for any
$(d_{1},...,d_{m})\in D(m;\mathcal{S})$. Given two simplicial infinitesimal
spaces $D(m;\mathcal{S})$ and $D(m^{\prime};\mathcal{S}^{\prime})$, a mapping
$\varphi=(\varphi_{1},...,\varphi_{m^{\prime}}):D(m;\mathcal{S})\rightarrow
D(m^{\prime};\mathcal{S}^{\prime})$ is called a \textit{monomial mapping} if
every $\varphi_{j}$ is a monomial in $d_{1},...,d_{m}$ with coefficient $1$.

Given an infinitesimal space $\mathbb{E}$, a mapping $\gamma$ from
$\mathbb{E}$ to $M$ is called an $\mathbb{E}$-\textit{microcube} on $M$. We
denote by $\mathbf{T}^{\mathbb{E}}(M)$ the totality of $\mathbb{E}$-microcubes
on $M$. Given $x\in M$,\ we denote by $\mathbf{T}_{x}^{\mathbb{E}}(M)$ the
totality of $\mathbb{E}$-microcubes$\ \gamma$ on $M$ with$\ \gamma
(0,...,0)=x$. A mapping $f:M\rightarrow M^{\prime}$ of microlinear spaces
naturally gives rise to a mapping $\gamma\in\mathbf{T}^{\mathbb{E}}(M)\mapsto
f\circ\gamma\in\mathbf{T}^{\mathbb{E}}(M^{\prime})$, which is denoted by
$f_{\ast}^{\mathbb{E}}$ or $f_{\ast}$. It is well known that the canonical
projection $\tau_{E}:\mathbf{T}^{D}(E)\rightarrow E$ is a vector bundle. Its
subbundle $\nu_{\pi}:\mathbf{V}(\pi)\rightarrow E$ with $\mathbf{V}%
(\pi)=\{t\in\mathbf{T}^{D}(E):\pi_{\ast}(t)=0\}$, called the \textit{vertical
bundle} of $\pi$, is also a vector bundle.

%change: The following has been embedded into a sloppypar environment 
%to prevent the overfull box
\begin{sloppypar}Given $\gamma\in\mathbf{T}^{D_{n}}(M)$ and a simplicial infinitesimal space
$D(m;\mathcal{S})$ with $\mathrm{\dim}D(m;\mathcal{S})\leq n$, the mapping
$(d_{1},...,d_{m})\in D(m;\mathcal{S})\mapsto\gamma(d_{1}+...+d_{m})\in M$ is
denoted by $\gamma_{D(m;\mathcal{S})}$. The following lemma should be obvious.
\end{sloppypar}

\begin{lemma}
\label{t0.2.1}Let $\gamma\in\mathbf{T}^{D_{n}}(M)$. Let $D(m;\mathcal{S}%
_{1})\subseteq D(m;\mathcal{S}_{2})$ (this is the case if $\mathcal{S}%
_{1}\supseteq\mathcal{S}_{2}$) with $n_{1}=\mathrm{\dim}D(m;\mathcal{S}%
_{1})\leq n_{2}=\mathrm{\dim}D(m;\mathcal{S}_{2})\leq n$. Let $n_{1}\leq
n_{3}\leq n$. Then we have
\[
\gamma_{D(m;\mathcal{S}_{2})}\mid_{D(m;\mathcal{S}_{1})}=(\gamma\mid
_{D_{n_{3}}})_{D(m;\mathcal{S}_{1})}%
\]

\end{lemma}

We denote by $\mathfrak{S}_{n}$ the symmetric group of the set $\{1,...,n\}$,
which is well known to be generated by $n-1$ transpositions $<i,i+1>$
exchanging $i$ and $i+1(1\leq i\leq n-1)$ while keeping the other elements
fixed. Given $\sigma\in\mathfrak{S}_{n}$ and $\gamma\in\mathbf{T}_{x}^{D^{n}%
}(M)$, we define $\Sigma_{\sigma}(\gamma)\in\mathbf{T}_{x}^{D^{n}}(M)$ to be
\[
\Sigma_{\sigma}(\gamma)(d_{1},...,d_{n})=\gamma(d_{\sigma(1)},...,d_{\sigma
(n)})
\]
for any $(d_{1},...,d_{n})\in D^{n}$. Given $\alpha\in\mathbb{R}$ and
$\gamma\in\mathbf{T}^{D^{n}}(M)$, we define $\alpha\underset{i}{\cdot}%
\gamma\in\mathbf{T}_{x}^{D^{n}}(M)\ (1\leq i\leq n)$ to be
\[
(\alpha\underset{i}{\cdot}\gamma)(d_{1},...,d_{n})=\gamma(d_{1},...,d_{i-1}%
,\alpha d_{i},d_{i+1},...,d_{n})
\]
for any $(d_{1},...,d_{n})\in D^{n}$. Given $\alpha\in\mathbb{R}$ and
$\gamma\in\mathbf{T}^{D_{n}}(M)$, we define $\alpha\gamma\in\mathbf{T}%
_{x}^{D_{n}}(M)\ (1\leq i\leq n)$ to be
\[
(\alpha\gamma)(d)=\gamma(\alpha d)
\]
for any $d\in D_{n}$. For any $\gamma\in\mathbf{T}_{x}^{D_{n}}(M)$ and any
$d\in D_{n}$, we define $\mathbf{i}_{d}(\gamma)\in\mathbf{T}_{x}^{D_{n+1}%
}(M)\ $to be
\[
\mathbf{i}_{d}(\gamma)(d^{\prime})=\gamma(dd^{\prime})
\]
for any $d^{\prime}\in D_{n+1}$.

\subsection{Affine Bundles}

A bundle $\pi:E\rightarrow M$ is called an \textit{affine bundle} over a
vector bundle $\xi:P\rightarrow M$ if $E_{x}$ is an affine space over the
vector space $P_{x}$ for any $x\in M$. Following the lines of our previous
paper \cite{n4}, we will establish a variant of Theorem 0.3.6 of that paper.

As in Lemma 0.3.2 of \cite{n4}, we have

\begin{lemma}
\label{t0.3.4}The diagram
%change: improved diagram
\[ \begin{CD}
D_{n}@>i>> D_{n+1}\\
@ViVV @VV{\psi _{D_{n+1}}}V\\
D_{n+1}@>>\varphi _{D_{n+1}}> D_{n+1}\oplus D
\end{CD} \]
is a quasi-colimit diagram, where $\varphi_{D_{n+1}}(d)=(d,0)$ and
$\psi_{D_{n+1}}(d)=(d,d^{n+1})$ for any $d\in D_{n+1}$.
\end{lemma}

\begin{remark}
We will write $\varphi_{D^{n+1}}:D^{n+1}\rightarrow D^{n+1}\oplus D$ and
$\psi_{D^{n+1}}:D^{n+1}\rightarrow D^{n+1}\oplus D$ for the mappings referred
to in Lemma 0.3.2 of \cite{n4}.
\end{remark}

Given two $D_{n+1}$-microcubes$\ \gamma_{+}$ and $\gamma_{-}$ on $M$ with
$\gamma_{+}|_{D_{n}}=\gamma_{-}|_{D_{n}}$, there exists a unique function
$f:D_{n+1}\oplus D\rightarrow M$ with $f\circ\psi_{D_{n+1}}=\gamma_{+}$ and
$f\circ\varphi_{D_{n+1}}=\gamma_{-}$. We define $(\gamma_{+}\dot{-}\gamma
_{-})\in\mathbf{T}^{D}(M)$ to be
\[
(\gamma_{+}\dot{-}\gamma_{-})(d)=f(0,d)
\]
for any $d\in D$. It is easy to see that
\begin{equation}
\alpha\gamma_{+}\dot{-}\alpha\gamma_{-}=\alpha^{n+1}(\gamma_{+}\dot{-}%
\gamma_{-}) \label{0.3.2}%
\end{equation}
for any $\alpha\in\mathbb{R}$ and any $\gamma_{\pm}\in\mathbf{T}^{D_{n+1}}(M)$
with $\gamma_{+}|_{D_{n}}=\gamma_{-}|_{D_{n}}$, while we know well that%
\begin{equation}
\alpha\underset{i}{\cdot}\gamma_{+}\dot{-}\alpha\underset{i}{\cdot}\gamma
_{-}=\alpha(\gamma_{+}\dot{-}\gamma_{-})\text{ \ \ }(1\leq i\leq n+1)
\label{0.3.3}%
\end{equation}
for any $\alpha\in\mathbb{R}$ and any $\gamma_{\pm}\in\mathbf{T}^{D^{n+1}}(M)$
with $\gamma_{+}|_{D(n+1)_{n}}=\gamma_{-}|_{D(n+1)_{n}}$. By the very
definition of $\dot{-}$, we have

\begin{proposition}
\label{t0.3.7}Let $f:M\rightarrow M^{\prime}$. Given $\gamma_{+},\gamma_{-}%
\in\mathbf{T}^{D_{n+1}}(M)$ with $\gamma_{+}|_{D_{n}}=\gamma_{-}|_{D_{n}}$, we
have $f_{\ast}(\gamma_{+})|_{D_{n}}=f_{\ast}(\gamma_{-})|_{D_{n}}\ $and%
\[
f_{\ast}(\gamma_{+}\dot{-}\gamma_{-})=f_{\ast}(\gamma_{+})\dot{-}f_{\ast
}(\gamma_{-})
\]

\end{proposition}

As in Lemma 0.3.4 of \cite{n4}, we have

%change: new diagram
\begin{lemma}
\label{t0.3.5}The diagram
\[ \begin{CD}
1@>i>> D\\
@ViVV @VV{\epsilon _{D_{n+1}}}V\\
D_{n+1}@>>{\varphi _{D_{n+1}}}> D_{n+1}\oplus D
\end{CD} \]
is a quasi-colimit diagram, where $\varepsilon_{D_{n+1}}(d)=(0,d)$ for any
$d\in D$.
\end{lemma}

\begin{remark}
We will write $\varepsilon_{D^{n+1}}:D\rightarrow D^{n+1}\oplus D$ for the
mapping referred to in Lemma 0.3.4 of \cite{n4}.
\end{remark}

Given $t\in\mathbf{T}^{D}(M)\ $and $\gamma\in\mathbf{T}^{D_{n+1}}(M)\ $with
$t(0)=\gamma(0)$, there exists a unique function $f:D_{n+1}\oplus D\rightarrow
M\ $with $f\circ\varphi_{D_{n+1}}=\gamma\ $and $f\circ\varepsilon_{D_{n+1}}%
=t$. We define $(t\dot{+}\gamma)\in\mathbf{T}^{D_{n+1}}(M)$ to be
\[
(t\dot{+}\gamma)(d)=f(d,d^{n+1})
\]
for any $d\in D_{n+1}$. It is easy to see that%
\begin{equation}
\alpha^{n+1}t\dot{+}\alpha\gamma=\alpha(t\dot{+}\gamma) \label{0.3.4}%
\end{equation}
for any $\alpha\in\mathbb{R}$. By the very definition of $\dot{+}$ we have

\begin{proposition}
\label{t0.3.8}Let $f:M\rightarrow M^{\prime}$. Given $t\in\mathbf{T}^{D}%
(M)\ $and $\gamma\in\mathbf{T}^{D_{n+1}}(\gamma)$with $t(0)=\gamma(0,...,0)$,
we have $f_{\ast}(t)(0)=f_{\ast}(\gamma)(0,...,0)\ $and%
\[
f_{\ast}(t\dot{+}\gamma)=f_{\ast}(t)\dot{+}f_{\ast}(\gamma)
\]

\end{proposition}

As in Theorem 0.3.6 of \cite{n4}, we have the following affine bundle theorem.

\begin{theorem}
\label{t0.3.6}The canonical projection $\mathbf{T}^{D_{n+1}}(M)\rightarrow
\mathbf{T}^{D_{n}}(M)\ $is an affine bundle over the vector bundle
$\mathbf{T}^{D}(M)\underset{M}{\times}\mathbf{T}^{D_{n}}(M)\rightarrow
\mathbf{T}^{D_{n}}(M)$.
\end{theorem}

\subsection{The First Approach to Jet Bundles}

In \cite{n4} we have discussed a synthetic approach to jet bundles in terms of
infinitesimal spaces $D^{n}$'s. We have no intention to reproduce the paper,
but our present notation and terminology are slightly different from those of
that paper. First of all, what were called $n$-pseudoconnections and
$n$-preconnections in that paper are now called $D^{n}$%
\textit{-pseudotangentials} and $D^{n}$\textit{-tangentials} respectively. We
will write $\mathbb{J}_{x}^{D^{n}}(\pi)$ and $\mathbb{J}^{D^{n}}(\pi)$ for the
space of $D^{n}$-tangentials for $\pi:E\rightarrow M$ at $x\in E$ and that of
$D^{n}$-tangentials for $\pi:E\rightarrow M$ in place of $\mathbb{J}_{x}%
^{n}(\pi)$ and $\mathbb{J}^{n}(\pi)$ in that paper respectively. We will
denote the canonical projection $\mathbb{J}^{D^{n+1}}(\pi)\rightarrow
\mathbb{J}^{D^{n}}(\pi)$ by $\pi_{n+1,n}$ in place of \underline{$\pi$%
}$_{n+1,n}$ in that paper. What is more than a matter of notation or
terminology is that, in our present definition of a $D^{n+1}$-tangential
$\mathsf{f}$ for $\pi:E\rightarrow M$ at $x\in E$, we require the following
condition besides conditions (1.11) and (1.12) of our previous paper
\cite{n4}:%
\begin{align*}
\mathsf{f}((d_{1},...,d_{n+1})  &  \in D^{n+1}\mapsto\gamma(ed_{1}%
,...,ed_{n+1})\in M)\\
&  =(d_{1},...,d_{n+1})\in D^{n+1}\mapsto\mathsf{f}_{D(n+1)_{n}}%
(\gamma)(ed_{1},...,ed_{n+1})\in E
\end{align*}
for any $e\in D_{n}$ and any $\gamma\in\mathbf{T}_{\pi(x)}^{D(n+1)_{n}}(M)$,
where $\mathsf{f}_{D(n+1)_{n}}$ stands for the induced $D(n+1)_{n}%
$-pseudotangential of $\mathsf{f}$ depicted in Theorem 7 of our previous paper
\cite{n5}.

\subsection{Symmetric Forms}

Given a vector bundle $\xi:P\rightarrow E$, a \textit{symmetric} $D_{n}%
$-\textit{form}$\mathit{\ }$\textit{at} $x\in E$ \textit{with} \textit{values}
\textit{in} $\xi\ $is a mapping $\omega:\mathbf{T}_{\pi(x)}^{D_{n}%
}(M)\rightarrow P_{x}\ $subject to the following conditions:

\begin{enumerate}
\item For any $\gamma\in\mathbf{T}_{\pi(x)}^{D_{n}}(M)$ and any $\alpha
\in\mathbb{R}$,$\ $we have
\begin{equation}
\omega(\alpha\gamma)=\alpha^{n}\omega(\gamma) \label{0.4.2.1}%
\end{equation}

\item For any simple polynomial $\rho$ of $d\in D_{n}$ and any $\gamma
\in\mathbf{T}_{\pi(x)}^{D_{l}}(M)$ with $\mathrm{\dim}_{n}\rho=l<n$, we have
\begin{equation}
\omega(\gamma\circ\rho)=0 \label{0.4.2.2}%
\end{equation}

\end{enumerate}

We denote by $\mathbb{S}_{x}^{D_{n}}(\pi;\xi)\ $the totality of symmetric
$D_{n}$-\textit{forms}$\ $at $x\ $with values in $\xi$. We denote by
$\mathbb{S}^{D_{n}}(\pi;\xi)\ $the\ set-theoretic union of $\mathbb{S}%
_{x}^{D_{n}}(\pi;\xi)$'s$\ $for all $x\in E$. The canonical projection
$\mathbb{S}^{^{D_{n}}}(\pi;\xi)\rightarrow E$ is obviously a vector bundle.

What were simply called symmetric $n$-forms along $\pi$ with values in $\xi$
in our previous paper \cite{n4} are now called symmetric $D^{n}$%
-\textit{forms}$\ $with values in $\xi$ in distinction to symmetric $D_{n}%
$-forms with values in $\xi$. Corresponding to our modification of the notion
of $D^{n+1}$-tangential in the preceding subsection, a $D^{n}$-form $\omega$
is required to satisfy%
\[
\omega((d_{1},...,d_{n})\in D^{n}\mapsto\gamma(ed_{1},...,ed_{n})\in M)=0
\]
for any $e\in D_{n-1}$ and any $\gamma\in\mathbf{T}_{\pi(x)}^{D(n)_{n-1}}(M)$
besides conditions (0.4.1)-(0.4.3) of our previous paper \cite{n4}. The spaces
$\mathbb{S}_{x}^{D^{n}}(\pi;\xi)$ and $\mathbb{S}^{D^{n}}(\pi;\xi)$ shall be
such as are expected.

\section{The Second Approach to Jet Bundles}

We have already discussed the second approach to jet bundles in terms of
infinitesimal spaces $D_{n}$'s\ in \cite{n5}. What has remained there to be
discussed is the affine bundle theorem as a $D_{n}$-variant of Theorem 1.9 in
\cite{n4}. To begin with, let us recall the fundamental concepts of our second approach.

\begin{definition}
\label{t10.3.1}A $D_{n}$\textrm{-pseudotangential for }$\pi$\textrm{\ at
}$x\in E$ is a mapping $\mathfrak{f}:\mathbf{T}_{\pi(x)}^{D_{n}}%
(M)\rightarrow\mathbf{T}_{x}^{D_{n}}(E)\ $subject to the following two conditions:

\begin{enumerate}
\item We have
\begin{equation}
\pi\circ(\mathfrak{f}(\gamma))=\gamma\label{10.3.1}%
\end{equation}
for $\gamma\in\mathbf{T}_{\pi(x)}^{D_{n}}(M)$.

\item We have
\begin{equation}
\mathfrak{f}(\alpha\gamma)=\alpha\mathfrak{f}(\gamma) \label{10.3.2}%
\end{equation}
for any $\gamma\in\mathbf{T}_{\pi(x)}^{D_{n}}(M)$ and any $\alpha\in
\mathbb{R}$.
\end{enumerate}
\end{definition}

We denote by $\mathbb{\hat{J}}_{x}^{D_{n}}(\pi)$ the totality of $D_{n}%
$-pseudotangentials for $\pi$ at $x\in E$. We denote by $\mathbb{\hat{J}%
}^{D_{n}}(\pi)$ the totality of $D_{n}$-pseudotangentials for $\pi$, i.e., the
set-theoretic union of $\mathbb{\hat{J}}_{x}^{D_{n}}(\pi)$'s for all $x\in E$.
We have the canonical projection $\hat{\pi}_{n+1,n}:\mathbb{\hat{J}}^{D_{n+1}%
}(\pi)\rightarrow\mathbb{\hat{J}}^{D_{n}}(\pi)$ such that%

\[
\mathfrak{f}(\mathbf{i}_{d}(\gamma))=\mathbf{i}_{d}(\hat{\pi}_{n+1,n}%
(\mathfrak{f})(\gamma))
\]
for any $\mathfrak{f}\in\mathbb{\hat{J}}^{D_{n+1}}(\pi)$, any $d\in D_{n}$ and
any $\gamma\in\mathbf{T}^{D_{n}}(M)$, for which the reader is referred to
Proposition 15 of our previous paper \cite{n5}. By assigning $\pi(x)\in M$ to
each $D^{n}$-pseudotangential for $\pi:E\rightarrow M$ at $x\in E$ we have the
canonical projection $\hat{\pi}_{n}:\mathbb{\hat{J}}^{D_{n}}(\pi)\rightarrow$
$M$, which is easily seen to be a vector bundle providing that $\pi$ is
already a vector bundle. Note that $\hat{\pi}_{n}\circ\hat{\pi}_{n+1,n}%
=\hat{\pi}_{n+1}$. For any natural numbers $n$, $m$ with $m\leq n$, we define
$\hat{\pi}_{n,m}:\mathbb{\hat{J}}^{D_{n}}(\pi)\rightarrow\mathbb{\hat{J}%
}^{D_{m}}(\pi)$ to be the composition $\hat{\pi}_{m+1,m}\circ...\circ\hat{\pi
}_{n,n-1}$.

\begin{definition}
\label{t10.3.2}The notion of a $D_{n}$\textrm{-tangential for }$\pi
$\textrm{\ at }$x\in E$ is defined inductively on $n$. The notion of a $D_{0}%
$-\textit{tangential} for $\pi$ at $x\in E$ and that of a $D_{1}%
$-\textit{tangential} for $\pi$ at $x\in E$ shall be identical with that of a
$D_{0}$-pseudotangential for $\pi$ at $x\in E$ and that of a $D_{1}%
$-pseudotangential for $\pi$ at $x\in E$ respectively. Now we proceed by
induction on $n$. A $D_{n+1}$-pseudotangential $\mathfrak{f}:\mathbf{T}%
_{\pi(x)}^{D_{n+1}}(M)\rightarrow\mathbf{T}_{x}^{D_{n+1}}(E)$ for $\pi$ at
$x\in E$ is called a $D_{n+1}$\textrm{-tangential for }$\pi$\textrm{\ at }$x$
if it acquiesces in the following two conditions:

\begin{enumerate}
\item $\hat{\pi}_{n+1,n}(\mathfrak{f})$ is a $D_{n}$-tangential for $\pi$ at
$x$.

\item For any simple polynomial $\rho$ of $d\in D_{n+1}$ with $\mathrm{\dim
}_{n}\rho=l\leq n$ and any $\gamma\in\mathbf{T}_{\pi(x)}^{D_{l}}(M)$, we have
\begin{align*}
&  \mathfrak{f}(\gamma\circ\rho)\\
&  =(\hat{\pi}_{n+1,l}(\mathfrak{f})(\gamma))\circ\rho
\end{align*}

\end{enumerate}
\end{definition}

We denote by $\mathbb{J}_{x}^{D_{n}}(\pi)$ the totality of $D_{n}$-tangentials
for $\pi$ at $x$. We denote by $\mathbb{J}_{x}^{D_{n}}(\pi)$ the totality of
$D_{n}$-tangentials for $\pi$, namely, the set-theoretic union of
$\mathbb{J}_{x}^{D_{n}}(\pi)$'s for all $x\in E$. By the very definition of a
$D_{n}$-tangential for $\pi$, the projection $\hat{\pi}_{n+1,n}:\mathbb{\hat
{J}}^{D_{n+1}}(\pi)\rightarrow\mathbb{\hat{J}}^{D_{n}}(\pi)$ is naturally
restricted to a mapping $\pi_{n+1,n}:\mathbb{J}^{D_{n+1}}(\pi)\rightarrow
\mathbb{J}^{D_{n}}(\pi)$. Similarly for $\pi_{n}:\mathbb{J}^{D_{n}}%
(\pi)\rightarrow M$ and $\pi_{n,m}:\mathbb{J}^{D_{n}}(\pi)\rightarrow
\mathbb{J}^{D_{m}}(\pi)$ with $m\leq n$. We can see easily that $\pi
_{n}:\mathbb{J}^{D_{n}}(\pi)\rightarrow M$ is naturally a vector bundle
providing that $\pi$ is already a vector bundle.

Now we will establish a $D_{n}$-variant of Theorem 1.9 of \cite{n4}. Let us
begin with a $D_{n}$-variant of Proposition 1.7 in \cite{n4}.

\begin{proposition}
\label{t2.3}Let $\mathfrak{f}^{+}$, $\mathfrak{f}^{-}\in\mathbb{J}%
_{x}^{D_{n+1}}(\pi)$ with $\pi_{n+1,n}(\mathfrak{f}^{+})=\ \pi_{n+1,n}%
(\mathfrak{f}^{-})$. Then the assignment $\gamma\in\mathbf{T}_{\pi
(x)}^{D_{n+1}}(M)\longmapsto\mathfrak{f}^{+}(\gamma)\dot{-}\mathfrak{f}%
^{-}(\gamma)$ belongs to $\mathbb{S}_{\pi(x)}^{D_{n+1}}(\pi;v_{\pi})$.
\end{proposition}

\begin{proof}
Since we have%
\begin{align*}
&  \pi_{\ast}(\mathfrak{f}^{+}(\gamma)\dot{-}\mathfrak{f}^{-}(\gamma))\\
&  =\pi_{\ast}(\mathfrak{f}^{+}(\gamma))\dot{-}\pi_{\ast}(\mathfrak{f}%
^{-}(\gamma))\\
&  \text{[By Proposition \ref{t0.3.7}]}\\
&  =0
\end{align*}
$\mathfrak{f}^{+}(\gamma)\dot{-}\mathfrak{f}^{-}(\gamma)$ belongs in
$\mathbf{V}_{x}^{1}(\pi)$. For any $\alpha\in\mathbb{R}$ we have%
\begin{align*}
&  \mathfrak{f}^{+}(\alpha\gamma)\dot{-}\mathfrak{f}^{-}(\alpha\gamma)\\
&  =\alpha\mathfrak{f}^{+}(\gamma)\dot{-}\alpha\mathfrak{f}^{-}(\gamma)\\
&  =\alpha^{n+1}(\mathfrak{f}^{+}(\gamma)\dot{-}\mathfrak{f}^{-}(\gamma))\\
&  \text{[By (\ref{0.3.2})]}%
\end{align*}
which implies that the assignment $\gamma\in\mathbf{T}_{\pi(x)}^{D_{n+1}%
}(M)\longmapsto\mathfrak{f}^{+}(\gamma)\dot{-}\mathfrak{f}^{-}(\gamma
)\in\mathbf{V}_{x}^{1}(\pi)\ $abides by (\ref{0.4.2.1}). It remains to show
that the assignment $\gamma\in\mathbf{T}_{\pi(x)}^{n+1}(M)\longmapsto
\mathfrak{f}^{+}(\gamma)\dot{-}\mathfrak{f}^{-}(\gamma)\in\mathbf{V}_{x}%
^{1}(\pi)$ abides by (\ref{0.4.2.2}), which follows directly from the second
condition in Definition \ref{t10.3.2} and the assumption that $\pi
_{n+1,n}(\mathfrak{f}^{+})=\ \pi_{n+1,n}(\mathfrak{f}^{-})$.
\end{proof}

Now we will establish a $D_{n}$-variant of Proposition 1.8 in \cite{n4}.

\begin{proposition}
\label{t2.4}Let $\mathfrak{f}\in\mathbb{J}_{x}^{D_{n+1}}(\pi)$ and $\omega
\in\mathbb{S}_{\pi(x)}^{D_{n+1}}(\pi;v_{\pi})$. Then the assignment $\gamma
\in\mathbf{T}_{\pi(x)}^{D_{n+1}}(M)\longmapsto\omega(\gamma)\dot
{+}\mathfrak{f}(\gamma)$ belongs to $\mathbb{J}_{x}^{D_{n+1}}(\pi)$.
\end{proposition}

\begin{proof}
Since we have%
\begin{align*}
&  \pi_{\ast}(\omega(\gamma)\dot{+}\mathfrak{f}(\gamma))\\
&  =\pi_{\ast}(\omega(\gamma))\dot{+}\pi_{\ast}(\mathfrak{f}(\gamma))\\
&  \text{[By Proposition \ref{t0.3.8}]}\\
&  =\gamma
\end{align*}
the assignment $\gamma\in\mathbf{T}_{\pi(x)}^{D_{n+1}}(M)\longmapsto
\omega(\gamma)\dot{+}\mathfrak{f}(\gamma)$ stands to (\ref{10.3.1}). For any
$\alpha\in\mathbb{R}$ we have%
\begin{align*}
&  \omega(\alpha\gamma)\dot{+}\mathfrak{f}(\alpha\gamma)\\
&  =\alpha^{n+1}\omega(\gamma)\dot{+}\alpha\mathfrak{f}(\gamma)\\
&  \text{[By (\ref{0.3.4})]}\\
&  =\alpha(\omega(\gamma)\dot{+}\mathfrak{f}(\gamma))
\end{align*}
so that the assignment $\gamma\in\mathbf{T}_{\pi(x)}^{D_{n+1}}(M)\longmapsto
\omega(\gamma)\dot{+}\mathfrak{f}(\gamma)$ stands to (\ref{10.3.2}). That the
assignment $\gamma\in\mathbf{T}_{\pi(x)}^{D_{n+1}}(M)\longmapsto\omega
(\gamma)\dot{+}\mathfrak{f}(\gamma)$ stands to the first condition of
Definition \ref{t10.3.2} follows from the simple fact that the image of the
assignment under $\hat{\pi}_{n+1,n}$ coincides with $\hat{\pi}_{n+1,n}%
(\mathfrak{f})$, which is consequent upon (\ref{0.4.2.2}). It remains to show
that the assignment $\gamma\in\mathbf{T}_{\pi(x)}^{D_{n+1}}(M)\longmapsto
\omega(\gamma)\dot{+}\mathfrak{f}(\gamma)$ abides by the second condition of
Definition \ref{t10.3.2}, which follows directly from (\ref{0.4.2.2}) and the
assumption that $\mathfrak{f}$ satisfies the second condition of Definition
\ref{t10.3.2}.
\end{proof}

For any $\mathfrak{f}^{+},\mathfrak{f}^{-}\in\mathbb{J}^{D_{n+1}}(\pi)$ with
$\pi_{n+1,n}(\mathfrak{f}^{+})=\pi_{n+1,n}(\mathfrak{f}^{-})$, we define
$\mathfrak{f}^{+}\dot{-}\mathfrak{f}^{-}\in\mathbb{S}^{D_{n+1}}(M;v_{\pi})$ to
be
\[
(\mathfrak{f}^{+}\dot{-}\mathfrak{f}^{-})(\gamma)=\mathfrak{f}^{+}(\gamma
)\dot{-}\mathfrak{f}^{-}(\gamma)
\]
for any $\gamma\in\mathbf{T}_{\pi(x)}^{D_{n+1}}(M)$. For any $\mathfrak{f}%
\in\mathbb{J}_{x}^{D_{n+1}}(\pi)\ $and any $\omega\in\mathbb{S}_{\pi
(x)}^{D_{n+1}}(\pi;v_{\pi})$ we define $\omega\dot{+}\mathfrak{f}\in
\mathbb{J}_{x}^{D_{n+1}}(\pi)$ to be
\[
(\omega\dot{+}\mathfrak{f})(\gamma)=\omega(\gamma)\dot{+}\mathfrak{f}(\gamma)
\]
for any $\gamma\in\mathbf{T}_{\pi(x)}^{D_{n+1}}(M)$.

With these two operations we have the following microlinear generalization of
the classical affine bundle theorem (cf. Theorem 6.2.9 of Saunders \cite{s1}),
which has been concerned merely with the finite-dimensional setting.

\begin{theorem}
\label{t2.5}The bundle $\pi_{n+1,n}:\mathbb{J}^{D_{n+1}}(\pi)\rightarrow
\mathbb{J}^{D_{n}}(\pi)$ is an affine bundle over the vector bundle
$\underset{}{\mathbb{S}^{D_{n+1}}(\pi;v_{\pi})\underset{E}{\times}}%
\mathbb{J}^{D_{n}}(\pi)\rightarrow\mathbb{J}^{D_{n}}(\pi)$.
\end{theorem}

\begin{proof}
This follows simply from Theorem \ref{t0.3.6}.
\end{proof}

\section{The Comparison without Coordinates}

The relationship between $\mathbb{J}^{D^{n}}(\pi)$ and $\mathbb{J}^{D_{n}}%
(\pi)$ without any reference to the affine bundle structures stated in Theorem
1.9 of our \cite{n4} and Theorem \ref{t2.5} of this paper has already been
discussed in our \cite{n5}. We remind the reader that all the results of our
previous paper \cite{n5} are valid with due modifications, though our
definition of a simple polynomial and that of a $D^{n}$-tangential have been modified.

Let $\mathsf{f}$ be a $D^{n}$-tangential for $\pi:E\rightarrow M$ at $x\in E$.
We have a function $\Phi_{n}(\mathsf{f}):\mathbf{T}_{\pi(x)}^{D_{n}%
}(M)\rightarrow\mathbf{T}_{x}^{D_{n}}(E)$, which is characterized by%
\[
\mathsf{f}(\gamma_{D^{n}})=\Phi_{n}(\mathsf{f})(\gamma)_{D^{n}}%
\]
for any $\gamma\in\mathbf{T}_{\pi(x)}^{D_{n}}(M)$. For the existence and
uniqueness of such a function, the reader is referred to Lemma 18 of our
previous \cite{n5}, from which we quote two crucial results.

\begin{theorem}
\label{t3.5}For any\textbf{\ }$\mathsf{f}\in\mathbb{J}_{x}^{D^{n}}(\pi)$, we
have $\Phi_{n}(\mathsf{f})\in\mathbb{J}_{x}^{D_{n}}(\pi)$, so that we have a
canonical mapping $\Phi_{n}:\mathbb{J}^{D^{n}}(\pi)\rightarrow\mathbb{J}%
^{D_{n}}(\pi)$.
\end{theorem}

\begin{proposition}
\label{t3.6}The diagram
% change: improved diagram
\[ \begin{CD}
{\mathbb J}_{x }^{D^{n+1}} (\pi ) @>{\Phi _{n+1} }>>{\mathbb J}_{x}^{D_{n+1}}(\pi)\\
@V\pi _{n+1,n} VV @VV{\pi _{n+1,n}}V\\
\mathbb{J}_{x}^{D^{n}}(\pi )@>>{\Phi _{n}}> \mathbb{J}_{x}^{D_{n}}(\pi ) \end{CD}
\]
is commutative.
\end{proposition}

Now we are in a position to investigate the relationship between the affine
bundles discussed in Theorem 1.9 of our \cite{n4} and Theorem \ref{t2.5} of
this paper. Let us begin with

\begin{lemma}
\label{t3.1}Let $\gamma^{\pm}\in\mathbf{T}_{x}^{D_{n+1}}(E)$ with $\gamma
^{+}\mid_{D_{n}}=\gamma^{-}\mid_{D_{n}}$. Then
\begin{equation}
\gamma_{D^{n+1}}^{+}\mid_{D(n+1)_{n}}=\gamma_{D^{n+1}}^{-}\mid_{D(n+1)_{n}%
}\text{,}%
\end{equation}
and we have
\begin{equation}
\gamma^{+}\dot{-}\gamma^{-}=\gamma_{D^{n+1}}^{+}\dot{-}\gamma_{D^{n+1}}^{-}%
\end{equation}

\end{lemma}

\begin{proof}
Since $D(n+1)_{n}\subseteq D^{n+1}$, the first statement follows from the
following simple calculation:
\begin{align}
\gamma_{D^{n+1}}^{+}  &  \mid_{D(n+1)_{n}}\nonumber\\
&  =(\gamma^{+}\mid_{D_{n}})_{D(n+1)_{n}}\text{ \ \ [by Lemma \ref{t0.2.1}%
]}\nonumber\\
&  =(\gamma^{-}\mid_{D_{n}})_{D(n+1)_{n}}\nonumber\\
&  =\gamma_{D^{n+1}}^{-}\mid_{D(n+1)_{n}}\text{ \ \ [by Lemma \ref{t0.2.1}]}%
\end{align}
The second statement follows simply from the following commutative diagram%
\[%
%TCIMACRO{\FRAME{itbpF}{2.8426in}{1.7538in}{0in}{}{}{fig3.jpg}%
%width 2.8426in;  height 1.7538in;  depth 0in;  original-width 3.5924in;
%original-height 2.2018in;  cropleft "0";  croptop "1";  cropright "1";
%filename '../Local Settings/Temporary Internet Files/Content.IE5/CPIFCH6R/Fig3.jpg';file-properties "NPEU";}%
%}}%
%BeginExpansion
{\includegraphics[
natheight=2.201800in,
natwidth=3.592400in,
height=1.7538in,
width=2.8426in
]%
{Fig3.pdf}%
}%
%EndExpansion
\]
where $+$\ stands for addition of components.
\end{proof}

\begin{lemma}
\label{t3.2}Let $t\in\mathbf{T}_{x}^{1}(E)$ and $\gamma\in\mathbf{T}%
_{x}^{D_{n}}(E)$. Then we have
\begin{equation}
(t\dot{+}\gamma)_{D^{n}}=t\dot{+}\gamma_{D^{n}}%
\end{equation}

\end{lemma}

\begin{proof}
This follows simply from the following commutative diagram.%
\[%
%TCIMACRO{\FRAME{itbpF}{2.7985in}{1.7538in}{0in}{}{}{fig4.jpg}%
%width 2.7985in;  height 1.7538in;  depth 0in;  original-width 3.531in;
%original-height 2.2001in;  cropleft "0";  croptop "1";  cropright "1";
%filename '../Local Settings/Temporary Internet Files/Content.IE5/CPIFCH6R/Fig4.jpg';file-properties "NPEU";}%
%}}%
%BeginExpansion
{\includegraphics[
natheight=2.200100in,
natwidth=3.531000in,
height=1.7538in,
width=2.7985in
]%
{Fig4.pdf}%
}%
%EndExpansion
\]

\end{proof}

Now we are ready to state the main result of this section.

\begin{theorem}
\label{t3.3}We have the following.

\begin{enumerate}
\item For any $\mathsf{f}^{+},\mathsf{f}^{-}\in\mathbf{J}_{x}^{n}(\pi)$ and
any $\gamma\in\mathbf{T}_{\pi(x)}^{D_{n}}(M)$, we have
\begin{equation}
\Phi_{n}(\mathsf{f}^{+})(\gamma)\dot{-}\Phi_{n}(\mathsf{f}^{-})(\gamma
)=\mathsf{f}^{+}(\gamma_{D^{n}})\dot{-}\mathsf{f}^{-}(\gamma_{D^{n}})
\end{equation}

\item For any $\mathsf{f}\in\mathbf{J}_{x}^{n}(\pi)$, any $t\in\mathbf{T}%
_{\pi(x)}^{D}(M)$ and any $\gamma\in\mathbf{T}_{\pi(x)}^{D_{n}}(M)$, we have
\begin{equation}
(t\dot{+}\Phi_{n}(\mathsf{f})(\gamma))_{D^{n}}=t\dot{+}\mathsf{f}%
(\gamma_{D^{n}})
\end{equation}

\end{enumerate}
\end{theorem}

\begin{proof}
\begin{enumerate}
\item Since $\mathsf{f}^{\pm}(\gamma_{D^{n}})=(\Phi_{n}(\mathsf{f}^{\pm
})(\gamma))_{D^{n}}$, we have
\begin{align}
&  \mathsf{f}^{+}(\gamma_{D^{n}})\dot{-}\mathsf{f}^{-}(\gamma_{D^{n}%
})\nonumber\\
&  =(\Phi_{n}(\mathsf{f}^{+})(\gamma))_{D^{n}}\dot{-}(\Phi_{n}(\mathsf{f}%
^{-})(\gamma))_{D^{n}}\nonumber\\
&  =\Phi_{n}(\mathsf{f}^{+})(\gamma)\dot{-}\Phi_{n}(\mathsf{f}^{-}%
)(\gamma)\text{ \ \ [by Lemma \ref{t3.1}]}%
\end{align}

\item Since $\mathsf{f}(\gamma_{D^{n}})=(\Phi_{n}(\mathsf{f})(\gamma))_{D^{n}%
}$, we have
\begin{align}
&  t\dot{+}\mathsf{f}(\gamma_{D^{n}})\nonumber\\
&  =t\dot{+}(\Phi_{n}(\mathsf{f})(\gamma))_{D^{n}}\nonumber\\
&  =(t\dot{+}\Phi_{n}(\mathsf{f})(\gamma))_{D^{n}}\text{ \ \ [by Lemma
\ref{t3.2}]}%
\end{align}

\end{enumerate}
\end{proof}

Now we would like to discuss the relationship between $\mathbb{S}^{D^{n}}%
(\pi;v_{\pi})$ and $\mathbb{S}^{D_{n}}(\pi;v_{\pi})$.

\begin{proposition}
\label{t3.7}For any $\omega\in\mathbb{S}_{x}^{D^{n}}(\pi;v_{\pi})$, the
mapping $\gamma\in\mathbf{T}_{x}^{D_{n}}(M)\mapsto\omega(\gamma_{D^{n}}) $,
denoted by $\Psi_{n}(\omega)$, belongs to $\mathbb{S}_{x}^{D_{n}}(\pi;v_{\pi
})$, thereby giving rise to a function $\Psi_{n}:\mathbb{S}^{D^{n}}(\pi
;v_{\pi})\rightarrow\mathbb{S}^{D_{n}}(\pi;v_{\pi})$.
\end{proposition}

\begin{proof}
For $n=0,1$, the statement is trivial. For any $\omega\in\mathbb{S}%
_{x}^{D^{n+1}}(M;v_{\pi})$, there exist $\mathsf{f}^{+},\mathsf{f}^{-}%
\in\mathbf{J}_{x}^{n+1}(\pi)$ such that $\pi_{n+1,n}(\mathsf{f}^{+}%
)=\pi_{n+1,n}(\mathsf{f}^{-})$ and $\omega=\mathsf{f}^{+}\dot{-}\mathsf{f}%
^{-}$. Then we have the following:

\begin{enumerate}
\item Let $\alpha\in\mathbb{R}$ and $\gamma\in\mathbf{T}_{x}^{D_{n+1}}(M) $.
Then we have
\begin{align*}
&  \omega((\alpha\gamma)_{D^{n+1}})\\
&  =\mathsf{f}^{+}((\alpha\gamma)_{D^{n+1}})\dot{-}\mathsf{f}^{-}%
((\alpha\gamma)_{D^{n+1}})\\
&  =\Phi_{n+1}(\mathsf{f}^{+})(\alpha\gamma)\dot{-}\Phi_{n+1}(\mathsf{f}%
^{-})(\alpha\gamma)\text{ \ \ \ [by Theorem \ref{t3.3}]}\\
&  =\alpha(\Phi_{n+1}(\mathsf{f}^{+})(\gamma))\dot{-}\alpha(\Phi
_{n+1}(\mathsf{f}^{-})(\gamma))\\
&  =\alpha^{n+1}(\Phi_{n+1}(\mathsf{f}^{+})(\gamma)\dot{-}\Phi_{n+1}%
(\mathsf{f}^{-})(\gamma))\\
&  =\alpha^{n+1}(\mathsf{f}^{+}(\gamma_{D^{n+1}})\dot{-}\mathsf{f}^{-}%
(\gamma_{D^{n+1}}))\text{ \ \ \ [by Theorem \ref{t3.3}]}\\
&  =\alpha^{n+1}\omega(\gamma_{D^{n+1}})
\end{align*}
so that $\Psi_{n+1}(\omega)$ abides by (\ref{0.4.2.1}).

\item Let $\rho$ be a simple polynomial of $d\in D_{n+1}$ and $\gamma
\in\mathbf{T}_{\pi(x)}^{D_{l}}(M)$ with $\mathrm{\dim}_{n+1}\rho=l<n+1$, we
have
\begin{align*}
&  \omega((\gamma\circ\rho)_{D^{n+1}})\\
&  =\mathsf{f}^{+}((\gamma\circ\rho)_{D^{n+1}})\dot{-}\mathsf{f}^{-}%
((\gamma\circ\rho)_{D^{n+1}})\\
&  =\Phi_{n+1}(\mathsf{f}^{+})(\gamma\circ\rho)\dot{-}\Phi_{n+1}%
(\mathsf{f}^{-})(\gamma\circ\rho)\text{ \ \ \ [by Theorem \ref{t3.3}]}\\
&  =(\pi_{n+1,l}(\Phi_{n+1}(\mathsf{f}^{+}))(\gamma))\circ\rho\dot{-}%
(\pi_{n+1,l}(\Phi_{n+1}(\mathsf{f}^{-}))(\gamma))\circ\rho\\
&  =(\Phi_{l}(\pi_{n+1,l}(\mathsf{f}^{+}))(\gamma))\circ\rho\dot{-}(\Phi
_{l}(\pi_{n+1,l}(\mathsf{f}^{-}))(\gamma))\circ\rho\text{ \ \ \ [by
Proposition \ref{t3.6}]}\\
&  =0\text{,}%
\end{align*}
so that $\Psi_{n+1}(\omega)$ abides by (\ref{0.4.2.2}).
\end{enumerate}
\end{proof}

Let us fix our terminology. Given an affine bundle $\pi_{1}:E_{1}\rightarrow
M_{1}$ over a vector bundle $\xi_{1}:P_{1}\rightarrow M_{1}$ and another
affine bundle $\pi_{2}:E_{2}\rightarrow M_{2}$ over another vector bundle
$\xi_{2}:P_{2}\rightarrow M_{2}$, a triple $(f,g,h)$ of mappings
$f:M_{1}\rightarrow M_{2}$, $g:E_{1}\rightarrow E_{2}$ and $h:P_{1}\rightarrow
P_{2}$ is called a \textit{morphism of affine bundles} from the affine bundle
$\pi_{1}:E_{1}\rightarrow M_{1}$ over the vector bundle $\xi_{1}%
:P_{1}\rightarrow M_{1}$ to the affine bundle $\pi_{2}:E_{2}\rightarrow E_{2}$
over the vector bundle $\xi_{2}:P_{2}\rightarrow M_{2}$ provided that they
satisfy the following conditions:

\begin{enumerate}
\item $(f,g)$ is a morphism of bundles from $\pi_{1}$ to $\pi_{2}$. In other
words, the following diagram is commutative:
\[%
\begin{array}
[c]{ccc}%
E_{1} & \underrightarrow{g} & E_{2}\\
\pi_{1}\downarrow &  & \downarrow\pi_{2}\\
M_{1} & \overrightarrow{f} & E_{2}%
\end{array}
\]

\item $(f,h)$ is a morphism of bundles from $\xi_{1}$ to $\xi_{2}$. In other
words, the following diagram is commutative:
\[%
\begin{array}
[c]{ccc}%
P_{1} & \underrightarrow{h} & P_{2}\\
\xi_{1}\downarrow &  & \downarrow\xi_{2}\\
M_{1} & \overrightarrow{f} & E_{2}%
\end{array}
\]

\item For any $x\in M_{1}$, $(g\mid_{E_{1,x}},h\mid_{P_{1,x}})$ is a morphism
of affine spaces from $(E_{1,x},P_{1,x})$ to $(E_{2,x},P_{2,x})$.
\end{enumerate}

Using this terminology, we can summarize Theorem \ref{t3.3} succinctly as follows:

\begin{theorem}
\label{t3.8}The triple $(\Phi_{n},\Phi_{n+1},\Psi_{n+1}\times\Phi_{n})$ of
mappings is a morphism of affine bundles from the affine bundle $\pi
_{n+1,n}:\mathbb{J}^{D^{n+1}}(\pi)\rightarrow\mathbb{J}^{D^{n}}(\pi)$ over the
vector bundle $\underset{}{\mathbb{S}^{D^{n+1}}(\pi;v_{\pi})\underset
{E}{\times}}\mathbb{J}^{D^{n}}(\pi)\rightarrow\mathbb{J}^{D^{n}}(\pi)$ in
Theorem 1.9 of \cite{n4} to the affine bundle $\pi_{n+1,n}:\mathbb{J}%
^{D_{n+1}}(\pi)\rightarrow\mathbb{J}^{D_{n}}(\pi)$ over the vector bundle
$\underset{}{\mathbb{S}^{D_{n+1}}(\pi;v_{\pi})\underset{E}{\times}}%
\mathbb{J}^{D_{n}}(\pi)\rightarrow\mathbb{J}^{D_{n}}(\pi)$ in Theorem
\ref{t2.5}.
\end{theorem}

\section{The Comparison with Coordinates}

Throughout this section we assume that the bundle $\pi:E\rightarrow M$ is a
formal bundle (cf. \cite{n3}), so that we can assume without any loss of
generality that $E=\mathbb{R}^{p+q}$, $M=\mathbb{R}^{p}$ and $\pi$ is the
canonical projection of $\mathbb{R}^{p+q}$\ to the first $p$ coordinates, for
our considerations to follow are always infinitesimal. We will let $i$ with or
without subscripts range over natural numbers between $1$ and $p$ (including
endpoints). By the general Kock axiom any $\gamma\in\mathbf{T}^{D^{n}}(M)$ is
of the form
\[
(d_{1},...,d_{n})\in D^{n}\longmapsto(x^{i})+\Sigma_{r=1}^{n}\Sigma_{1\leq
k_{1}<...<k_{r}\leq n}d_{k_{1}}...d_{k_{r}}(a_{k_{1}...k_{r}}^{i}%
)\in\mathbb{R}^{p}\text{, }%
\]
while any $\gamma\in\mathbf{T}^{D_{n}}(M)$ is of the form
\[
d\in D_{n}\longmapsto(x^{i})+\Sigma_{r=1}^{n}d^{r}(a_{\mathbf{r}}^{i}%
)\in\mathbb{R}^{p}\text{.}%
\]

The principal objective in this section is to show that

\begin{theorem}
\label{t4.1}For any natural number $n$, $\Phi_{n}:\mathbb{J}^{D^{n}}%
(\pi)\rightarrow\mathbb{J}^{D_{n}}(\pi)$ and $\Psi_{n}:\mathbb{S}^{D^{n}}%
(\pi;v_{\pi})\rightarrow\mathbb{S}^{D_{n}}(\pi;v_{\pi})$ are bijective.
\end{theorem}

We proceed by induction on $n$. For $n=0,1$, the theorem holds trivially. We
have shown in \cite{n5} that $\Phi_{n}$ is injective for any natural number
$n$, so that $\Psi_{n}$ is injective for any natural number $n$. With due
regard to Theorem \ref{t3.8}, it suffices to show that $\Psi_{n}$ is bijective
for any natural number $n$, for which we can and will use dimension counting
techniques. Let us remark the following two plain propositions, which may
belong to the folklore.

\begin{proposition}
\label{t4.2}Let $x=(x^{i})\in M$ and $V$ a finite-dimensional $\mathbb{R}%
$-module. Let $\omega:\mathbf{T}_{x}^{D^{n}}(M)\rightarrow V$ be a function
acquiescent in conditions (0.4.1) and (0.4.2) of our previous \cite{n4}. Then
$\omega$ is of the following form:
\begin{align*}
\omega((d_{1},...,d_{n})  &  \in D^{n}\longmapsto(x^{i})+\Sigma_{r=1}%
^{n}\Sigma_{1\leq k_{1}<...<k_{r}\leq n}d_{k_{1}}...d_{k_{r}}(a_{k_{1}%
...k_{r}}^{i})\in\mathbb{R}^{p})\\
&  =\Sigma\Omega_{\mathbf{J}_{1},...,\mathbf{J}_{s}}^{n}((a_{\mathbf{J}_{1}%
}^{i_{\mathbf{J}_{1}}}),...,(a_{\mathbf{j}_{_{s}}}^{i_{\mathbf{J}_{s}}%
}))\text{, }%
\end{align*}
where the last $\Sigma\ $is taken over all partitions of the
set$\ \{1,...,n\}\ $into nonempty subsets $\{\mathbf{J}_{1},...,\mathbf{J}%
_{s}\}$, $\Omega_{\mathbf{J}_{1},...,\mathbf{J}_{s}}^{n}:(\mathbb{R}^{p}%
)^{s}\rightarrow V$ is a symmetric $s$-linear mapping, and $a_{\mathbf{J}%
}^{i_{\mathbf{J}}}$ denotes $a_{k_{1}...k_{r}}^{i_{\mathbf{J}}}$ for
$\mathbf{J}=\{k_{1},...,k_{r}\}$ with $k_{1}<...<k_{r}$.
\end{proposition}

\begin{proof}
By the same token as in the proof of Proposition 11 (Section 1.2) of \cite{l1}.
\end{proof}

\begin{proposition}
\label{t4.3}Let $x\in M$ and $V$ a finite-dimensional $\mathbb{R}$-module. Let
$\omega:\mathbf{T}_{x}^{D_{n}}(M)\rightarrow V$ be a function acquiescent in
condition (\ref{0.4.2.1}). Then $\omega$ is of the following form:
\begin{align*}
\omega(d  &  \in D_{n}\longmapsto(x^{i})+\Sigma_{r=1}^{n}d^{r}(a_{\mathbf{r}%
}^{i})\in\mathbb{R}^{p})\\
&  =\Sigma(\Omega_{r_{1},...,r_{k}}^{n}((a_{\mathbf{r}_{1}}^{i_{1}%
}),...,(a_{\mathbf{r}_{k}}^{i_{k}})))\text{, }%
\end{align*}
where $\Omega_{r_{1},...,r_{k}}^{n}:(\mathbb{R}^{p})^{k}\rightarrow V$ is a
symmetric $k$-linear mapping, and the last $\Sigma\ $is taken over all
partitions of the positive integer $n\ $into positive integers $r_{1}%
,...,r_{k}$ (so that $n=r_{1}+...+r_{k}$) with $r_{1}\leq...\leq r_{k}$.
\end{proposition}

\begin{proof}
By the same token as in the proof of Proposition 11 (Section 1.2) of \cite{l1}.
\end{proof}

\begin{proposition}
\label{t4.4}For any $x\in M$, $\mathbb{R}$-modules $\mathbb{S}_{x}^{D^{n}}%
(\pi;v_{\pi})$ and $\mathbb{S}_{x}^{D_{n}}(\pi;v_{\pi})$ are of the same
dimension $q(
\begin{array}
[c]{c}%
p+n-1\\
n
\end{array}
)$, so that $\Psi_{n}:\mathbb{S}^{D^{n}}(\pi;v_{\pi})\rightarrow
\mathbb{S}^{D_{n}}(\pi;v_{\pi})$ is bijective.
\end{proposition}

\begin{proof}
Taking into consideration the condition (0.4.3) of our previous \cite{n4} in
Proposition \ref{t4.2}, we can conclude that the $\mathbb{R}$-module
$\mathbb{S}_{x}^{D^{n}}(\pi;v_{\pi})$ is of dimension $q(%
\begin{array}
[c]{c}%
p+n-1\\
n
\end{array}
)$, for $\Omega_{\mathbf{J}_{1},...,\mathbf{J}_{s}}^{n}$ is zero unless
$\{\mathbf{J}_{1},...,\mathbf{J}_{s}\}=\{\{1\},...,\{n\}\}$. Similarly, taking
into consideration the condition (\ref{0.4.2.2}) in Proposition \ref{t4.3}, we
can conclude that the $\mathbb{R}$-module $\mathbb{S}^{D_{n}}(\pi;v_{\pi})$ is
of dimension $q(%
\begin{array}
[c]{c}%
p+n-1\\
n
\end{array}
)$, for $\Omega_{r_{1},...,r_{k}}^{n}$ is zero except $\Omega_{1,...,1}^{n}$.
\end{proof}

\end{document}